\documentclass[10pt]{article}

\usepackage{amssymb}

\begin{document}

\author{David Carf\`{i}}
\title{Financial Lie groups}
\date{}
\maketitle

\begin{abstract}
In this paper we see the evolution of a capitalized financial event $e$,
with respect to a capitalization factor $f$, as the exponential map of a
suitably defined Lie group $G_{f,e}$, supported by the halfspace of
capitalized financial events having the same sign of $e$. The Lie group $%
G_{f,e}$ depends on the capitalization factor $f$ and on the event $e$
itself. After the extention of the definition of exponential map of a Lie
group, we shall eliminate the dependence on the financial event $e$,
recognizing the precence of essentialy one unique financial Lie semigroup,
supported by the entire space of capitalized financial events, determined by
the capitalization factor $f$.
\end{abstract}

\bigskip

\bigskip

\section{\textbf{Financial preliminaries}}

\bigskip

\textbf{Definition (of capitalized financial event).} \emph{We call \textbf{%
capitalized financial event} any triple }$e=(t,h,c)$\emph{\ of real numbers.
We call \textbf{zero financial event} any event of the type }$(t,h,0)$\emph{%
. Moreover we call \textbf{credits} (resp. \textbf{strict credits}) the
financial events with non negative (resp. positive) third component and 
\textbf{debts} (resp. \textbf{strict debts}) those with non positive (resp.
negative) third component. The event }$(t,h,-c)$ \emph{is called \textbf{the
opposite of} }$e=(t,h,c)$\emph{.}

\bigskip

\textbf{Financial interpretation.} We interpret $e=(t,h,c)$ as a financial
object characterized by

1) a \emph{reference time} $t$;

2) a \emph{capitalization time} $h$, meaning that the event $e$ exists and
is under capitalization (is in the financial market) since the time $t-h$
(we call this time \emph{origin of the capitalized event} $e$);

3) a \emph{capital} $c$, considered as the value of the event at the
reference time $t$.

\bigskip

\textbf{Definition.}\emph{\ If }$e=(t,h,c)$\emph{\ is a capitalized event,
we call:}

\emph{1) \textbf{reference time of the event} the real number }$t$\emph{;}

\emph{2) \textbf{capitalization time of} }$e$\emph{\ the real number }$h$%
\emph{;}

\emph{3) \textbf{capital (value) of} }$e$\emph{\ (at }$t$\emph{)\ the real }$%
c$\emph{;}

\emph{4) \textbf{time origin of} }$e$\emph{\ the time }$t-h$\emph{;}

\emph{5) \textbf{state of the event} }$e$\emph{\ the pair }$(h,c)$\emph{.}

\bigskip 

\textbf{Definition (fibrations of capitalized financial events).}\emph{\ The 
\textbf{state fibration of capitalized events} is the space }$\Bbb{R}^{3}$%
\emph{\ endowed with the first canonical projection }$pr_{1}$\emph{. In
other terms, we define the state fibration as the trivial fibration }$(\Bbb{R%
}^{3},\Bbb{R},pr_{1})$\emph{\ having as basis the affine space of reference
times. This trivial fibration is a vector bundle of type }$\Bbb{R}^{2}$\emph{%
\ and, for every time }$t$\emph{, the fiber }$pr_{1}^{-}(t)$\emph{\ is
called \textbf{the vector space (fiber) of capitalized financial states}.
Moreover we define \textbf{capital fibration of capitalized events} the
trivial fibration }$(\Bbb{R}^{3},\Bbb{R}^{2},pr_{12})$\emph{\ having as
basis the affine multi-time plane. This trivial fibration is a vector bundle
of type }$\Bbb{R}$\emph{\ and, for every pair }$(t,h)$\emph{, the fiber }$%
pr_{12}^{-}(t,h)$\emph{\ is called \textbf{the vector space (fiber) of
capitals}.}

\bigskip

\textbf{Definition (of global capitalization factor).} \emph{A (global) }$%
C^{1}$\emph{\ \textbf{capitalization factor} }$f$\emph{\ is any real
function defined over the real line of time displacements enjoing the
following properties:}

\emph{1) }$f$\emph{\ is positive;}

\emph{2) }$f$\emph{\ map the zero displacement }$0$\emph{\ into }$1$\emph{;}

\emph{3) }$f(-h)=f(h)^{-1}$\emph{, for every time displacement }$h$\emph{;}

\emph{4) }$f$\emph{\ is of class }$C^{1}$\emph{.}

\bigskip 

\emph{The pair }$(\Bbb{R}^{3},f)$\emph{\ is called \textbf{financial space
with capitalization factor} }$f$\emph{.}

\bigskip 

\section{\textbf{Lie product induced by a capitalization factor}}

\bigskip

We define the following algebraic product on the space of capitalized
financial events.

\bigskip

\textbf{Definition (Lie product induced by a capitalization factor).} \emph{%
Let }$f$\emph{\ be a capitalization factor over the displacement time line.
Let }$e=(t,h,c)$\emph{\ and }$e^{\prime }=(t^{\prime },h^{\prime },c^{\prime
})$\emph{\ be two capitalized financial events, we define their }$f$\emph{-%
\textbf{product} }$ee^{\prime }$\emph{\ to be the capitalized financial event%
} 
\[
ee^{\prime }=[e|e^{\prime }]_{f}=(t+t^{\prime },h+h^{\prime
},cf(-h)c^{\prime }f(-h^{\prime })f(h+h^{\prime })).
\]
\emph{Analogously, let }$s=(h,c)$\emph{\ and }$s^{\prime }=(h^{\prime
},c^{\prime })$\emph{\ be two capitalized financial states, we define their }%
$f$\emph{-\textbf{product} }$ss^{\prime }$\emph{\ to be the capitalized
financial state} 
\[
ss^{\prime }=[s|s^{\prime }]_{f}=(h+h^{\prime },cf(-h)c^{\prime
}f(-h^{\prime })f(h+h^{\prime })).
\]
\emph{We call the products so defined \textbf{Lie products induced by the
capitalization factor} }$f$\emph{.}

\bigskip 

We will use the brief notation $ee^{\prime }$ (or $ss^{\prime }$), instead
of the more precise $(ee^{\prime })_{f}$ or $[e|e^{\prime }]_{f}$, when no
confusion is possible.

\bigskip

\textbf{Remark 1.} The definition seems to be lacking from the dimensional
point of view and for what concerns the addition of two times (it is indeed
non proper to add two times) but we must regard the sum $t+t^{\prime }$ as
the sum of the time $t$ with the time duration $h=t^{\prime }-0$ (and this
is perfectly possible since the time line is a good affine space over the
real line of time displacements); moreover we have to consider the product $%
cc^{\prime }$ as the product of two capitals divided by $1$ monetary unity,
in order to obtain a third component of the product with the dimensions of a
capital.

\bigskip 

\section{\textbf{Lie anti-product induced by a capitalization factor}}

\bigskip

\textbf{Definition (Lie anti-product induced by a capitalization factor).} 
\emph{Let }$f$\emph{\ be a capitalization factor over the displacement time
line. Let }$e=(t,h,c)$\emph{\ and }$e^{\prime }=(t^{\prime },h^{\prime
},c^{\prime })$\emph{\ be two capitalized financial events, we define their }%
$f$\emph{-\textbf{anti-product} }$[e|e^{\prime }]_{(f,-)}$\emph{\ to be the
capitalized financial event} 
\[
\lbrack e|e^{\prime }]_{(f,-)}=(t+t^{\prime },h+h^{\prime },-cf(-h)c^{\prime
}f(-h^{\prime })f(h+h^{\prime })).
\]
\emph{Analogously, let }$s=(h,c)$\emph{\ and }$s^{\prime }=(h^{\prime
},c^{\prime })$\emph{\ be two capitalized financial states, we define their }%
$f$\emph{-\textbf{anti-product} }$[s|s^{\prime }]_{(f,-)}$\emph{\ to be the
capitalized financial state} 
\[
\lbrack s|s^{\prime }]_{(f,-)}=(h+h^{\prime },-cf(-h)c^{\prime }f(-h^{\prime
})f(h+h^{\prime })).
\]
\emph{We call the products so defined \textbf{Lie anti-products induced by
the capitalization factor} }$f$\emph{.}

\bigskip 

\textbf{Notation.} We shall denote the space $\Bbb{R}^{3}$ endowed by the
Lie product induced by the capitalization factor $f$, that is the structure $%
(\Bbb{R}^{3},[.|.]_{f})$, by the symbol $\Bbb{R}_{f}^{3}$ and we shall
denote the space $\Bbb{R}^{3}$ endowed by the Lie anti-product induced by
the capitalization factor $f$, that is the structure $(\Bbb{R}%
^{3},[.|.]_{(f,-)})$, by the symbol $\Bbb{R}_{(f,-)}^{3}$.

\bigskip 

\section{\textbf{Financial Lie groups}}

\bigskip

\textbf{Theorem 1.}\emph{\ Let }$f$\emph{\ be any }$C^{1}$\emph{\
capitalization factor. Then:}

\begin{itemize}
\item  \emph{the space }$\Bbb{R}^{3}$\emph{\ of capitalized events is a }$%
C^{1}$\emph{\ Abelian Lie semigroup, with respect to the standard Euclidean
differentiable structure on }$\Bbb{R}^{3}$\emph{\ and to the Lie product }$%
[.|.]_{f}$\emph{\ induced by the capitalization factor }$f$\emph{;}

\item  \emph{the orgin of the Lie semigroup }$\Bbb{R}_{f}^{3}$\emph{\ is the
point }$o=(0,0,1)$\emph{;}

\item  \emph{the subset of non zero capitalized event is a Lie subgroup of
the Lie semigroup }$\Bbb{R}_{f}^{3}$\emph{\ and coincides with the subgroup
of the invertible elements of the Lie semigroup }$\Bbb{R}_{f}^{3}$\emph{;}

\item  \emph{denoted the above subgroup, when endowed with the product }$%
[.|.]_{f}$\emph{, by }$G_{f}$\emph{; the Lie group }$G_{f}$\emph{\ has
exactly two connected homeomorphic components, one (the half space of
credits }$C$\emph{) contains the origin and is a subgroup of }$G_{f}$ \emph{%
and the other one (the half space of debts }$D$\emph{) contains the opposite
event }$(0,0,-1)$ \emph{of the origin of the Lie semigroup }$\Bbb{R}_{f}^{3}$
\emph{and is a group with respect to the anti-product induced by }$f$\emph{;}

\item  \emph{denoted by }$G_{f}^{>}$\emph{\ and }$G_{f}^{<}$\emph{\ the
above groups of credits and debts, respectively (pay attention, only }$%
G_{f}^{>}$\emph{\ is a subgroup of }$G_{f}$\emph{) the mapping sending }$%
(t,h,c)$\emph{\ into }$(t,h,-c)$\emph{\ is a Lie groups isomorphism of }$%
G_{f}^{>}$\emph{\ onto }$G_{f}^{<}$\emph{.}
\end{itemize}

\emph{\bigskip }

\emph{Proof.} We prove first that the product $[.|.]_{f}$ is a semigroup
product.

1) \emph{Associativity.} Let $e$, $e^{\prime }$ and $e^{\prime \prime }$ be
the events $(t,h,c)$, $(t^{\prime },h^{\prime },c^{\prime })$ and $%
(t^{\prime \prime },h^{\prime \prime },c^{\prime \prime })$, respectively.
Then we obtain (using the brief notation for the Lie product $[.|.]_{f}$) 
\begin{eqnarray*}
(ee^{\prime })e^{\prime \prime } &=&(t+t^{\prime },h+h^{\prime
},cf(-h)c^{\prime }f(-h^{\prime })f(h+h^{\prime }))e^{\prime \prime }= \\
&=&\left( t+t^{\prime }+t^{\prime \prime },h+h^{\prime }+h^{\prime \prime },%
\frac{c}{f(h)}\frac{c^{\prime }}{f(h^{\prime })}\frac{f(h+h^{\prime })}{%
f(h+h^{\prime })}\frac{c^{\prime \prime }}{f(h^{\prime \prime })}%
f(h+h^{\prime }+h^{\prime \prime })\right) = \\
&=&\left( t+t^{\prime }+t^{\prime \prime },h+h^{\prime }+h^{\prime \prime },%
\frac{c}{f(h)}\frac{c^{\prime }}{f(h^{\prime })}\frac{c^{\prime \prime }}{%
f(h^{\prime \prime })}\frac{f(h^{\prime }+h^{\prime \prime })}{f(h^{\prime
}+h^{\prime \prime })}f(h+h^{\prime }+h^{\prime \prime })\right) = \\
&=&e(t^{\prime }+t^{\prime \prime },h^{\prime }+h^{\prime \prime },c^{\prime
}f(-h^{\prime })c^{\prime \prime }f(-h^{\prime \prime })f(h^{\prime
}+h^{\prime \prime }))= \\
&=&e(e^{\prime }e^{\prime \prime }),
\end{eqnarray*}
as we desired. Note, by the way, that 
\[
ee^{\prime }e^{\prime \prime }=\left( t+t^{\prime }+t^{\prime \prime
},h+h^{\prime }+h^{\prime \prime },\frac{cc^{\prime }c^{\prime \prime
}f(h+h^{\prime }+h^{\prime \prime })}{f(h)f(h^{\prime })f(h^{\prime \prime })%
}\right) .
\]

2) \emph{Existence of neutral element.} The neutral element of the magma $%
\Bbb{R}_{f}^{3}$ is (obviously) the capitalized event $(0,0,1)$.

3) \emph{Commutativity.} The product $[.|.]_{f}$ is evidently commutative.

4) \emph{Continuous differentiability of the Lie product.} The Lie product
induced by $f$ is $C^{1}$ because $f$ is $C^{1}$. By the way, note that the
Lie product $[.|.]_{f}$, defined by 
\[
\lbrack e|e^{\prime }]_{f}=(t+t^{\prime },h+h^{\prime },cf(-h)c^{\prime
}f(-h^{\prime })f(h+h^{\prime })),
\]
is differentiable since each of his component is differentiable with respect
to each argument; for example, we have 
\[
\partial _{2}([.|.]_{f})(e,e^{\prime })=(0,1,-cf^{\prime }(-h)c^{\prime
}f(-h^{\prime })f^{\prime }(h+h^{\prime })),
\]
and 
\[
\partial _{3}([.|.]_{f})(e,e^{\prime })=(0,0,f(-h)c^{\prime }f(-h^{\prime
})f(h+h^{\prime })).
\]
5) \emph{Continuous differentiability of the Lie inverse map.} The
invertible events are all those events $e$ with non zero capital; indeed,
the inverse of a non zero event $(t,h,c)$ is the capitalized event $%
(-t,-h,1/c)$. Indeed, 
\[
(t,h,c)(-t,-h,1/c)=(0,0,cf(-h)c^{-1}f(h)f(h-h))=(0,0,1),
\]
note that the inverse $e_{f}^{-1}$ is independent of the capitalization
factor $f$. It is simple to prove that if $e$ is $f$-invertible then $%
pr_{3}(e)$ is different from $0$ (otherwise $0=1$). The inverse map 
\[
\lbrack .]^{-1}:\Bbb{R}^{3}\backslash N\rightarrow \Bbb{R}^{3}\backslash
N:e\mapsto e^{-1},
\]
where $N$ is the set of zero events (additive kernel of the third
projection), is a group isomorphism (independent of the capitalization
factor $f$) and its derivatives with respect to the arguments are 
\[
\partial _{1}([.]^{-1})(e)=(-1,0,0),\;\;\partial _{2}([.]^{-1})(e)=(0,-1,0),%
\;\;\partial _{3}([.]^{-1})(e)=(0,0,-1/c^{2}).
\]
6) The two connected components of $G_{f}$ are the sets of strict credits
and the set of stricts debts, indeed this two sets are obviously two
subgroups of $G$ and they are isomorphic by means of the opposite mapping
\[
(t,h,c)\mapsto (t,h,-c).
\]
The theorem is proved. $\blacksquare $

\bigskip 

We have the anti-version of the preceding theorem.

\bigskip 

\textbf{Theorem 2.}\emph{\ Let }$f$\emph{\ be any }$C^{1}$\emph{\
capitalization factor. Then:}

\begin{itemize}
\item  \emph{the space }$\Bbb{R}^{3}$\emph{\ of capitalized events is a }$%
C^{1}$\emph{\ Abelian Lie semigroup, with respect to the standard Euclidean
differentiable structure on }$\Bbb{R}^{3}$\emph{\ and to the Lie
anti-product }$[.|.]_{(f,-)}$\emph{\ induced by the capitalization factor }$f
$\emph{;}

\item  \emph{the orgin of the Lie semigroup }$\Bbb{R}_{f,-}^{3}$\emph{\ is
the point }$-o=(0,0,-1)$\emph{;}

\item  \emph{the subset of non zero capitalized event is a Lie subgroup of
the Lie semigroup }$\Bbb{R}_{f,-}^{3}$\emph{\ and coincides with the
subgroup of the invertible elements of the Lie semigroup }$\Bbb{R}_{f,-}^{3}$%
\emph{;}

\item  \emph{denoted the above subgroup, when endowed with the product }$%
[.|.]_{f,-}$\emph{, by }$G_{f,-}$\emph{; the Lie group }$G_{f,-}$\emph{\ has
exactly two connected homeomorphic components, one (the half space of debts }%
$D$\emph{) contains the origin }$-o$ \emph{and is a subgroup of }$G_{f,-}$ 
\emph{and the other one (the half space of credits }$C$\emph{) contains the
event }$o=(0,0,1)$\emph{, origin of the Lie semigroup }$\Bbb{R}_{f}^{3}$, 
\emph{and it is a group with respect to the product induced by }$f$\emph{;}

\item  \emph{denoted by }$G_{f}^{>}$\emph{\ and }$G_{f}^{<}$\emph{\ the
above groups of credits and debts, respectively (pay attention, only }$%
G_{f}^{<}$\emph{\ is a subgroup of }$G_{f,-}$\emph{) the mapping sending }$%
(t,h,c)$\emph{\ into }$(t,h,-c)$\emph{\ is a Lie groups isomorphism of }$%
G_{f}^{<}$\emph{\ onto }$G_{f}^{>}$\emph{.}
\end{itemize}

\bigskip 

\section{\textbf{The evolution of the unit event}}

\bigskip

Consider now the evolution of an event $e_{0}=(t_{0},h_{0},c_{0})$ with
respect to a capitalization factor $f$, that is (by definition) the curve 
\[
\mu _{e_{0}}:\Bbb{R}\rightarrow \Bbb{R}^{3}:t\mapsto
(t,h+t-t_{0},c_{0}f(-h_{0})f(h_{0}+t-t_{0})). 
\]

\bigskip

In particular, if the event $e_{0}$ is the unit event $o=(0,0,1)$, we have
simply 
\[
\mu _{e_{0}}:\Bbb{R}\rightarrow \Bbb{R}^{3}:t\mapsto (t,t,f(t)). 
\]

\bigskip

Let us see the first resul of the paper on financial dynamical systems.

\bigskip

\textbf{Theorem 3.}\emph{\ The evolution of the origin }$o$\emph{\ of the
semigroup }$\Bbb{R}_{f}^{3}$\emph{\ is the exponential map of }$\Bbb{R}%
_{f}^{3}$\emph{\ with respect to the tangent vector }$(o,(1,1,f^{\prime
}(0)))$\emph{, a tangent vector to the Lie semigroup }$\Bbb{R}_{f}^{3}$\emph{%
\ at the origin }$o$ \emph{itself.}

\emph{\bigskip }

\emph{Proof.} Note that $\mu _{o}$ is a one parameter group in the Lie group 
$G_{f}^{>}$ of credits (events with positive capital) with respect to the
Lie semigroup operation induced by the capitalization factor $f$. Indeed, we
have 
\[
\mu _{o}(t+t^{\prime })=(t+t^{\prime },t+t^{\prime },f(t+t^{\prime })), 
\]
for every pair $(t,t^{\prime })$ of times, and 
\begin{eqnarray*}
\lbrack \mu _{o}(t)|\mu _{o}(t^{\prime })]_{f} &=&(t,t,f(t))(t^{\prime
},t^{\prime },f(t^{\prime }))= \\
&=&(t+t^{\prime },t+t^{\prime },f(t)f(-t)f(t^{\prime })f(-t^{\prime
})f(t+t^{\prime }))= \\
&=&(t+t^{\prime },t+t^{\prime },f(t+t^{\prime })),
\end{eqnarray*}
again, for every $t$ and $t^{\prime }$ over the time line. Moreover, the
tangent vector at $t$ of the curve $\mu _{o}$ is 
\[
\mu _{o}^{\prime }(t)=(1,1,f^{\prime }(t)), 
\]
so that we have 
\[
\mu _{o}^{\prime }(0)=(1,1,f^{\prime }(0))=(1,1,\delta _{f}(0)), 
\]
where 
\[
\delta _{f}:\Bbb{R}\rightarrow \Bbb{R}:h\mapsto f^{\prime }(h)/f(h), 
\]
is the so called force of interest of the capitalization factor $f$. Now, as
it is well known in Lie Group Theory, there is only one $1$-parameter group
in a Lie group $G$ having a fixed tangent vector $v\in T_{o}(G)$ as a
tangent vector at $0$, and this is the exponential map 
\[
\exp _{v}:\Bbb{R}\rightarrow G, 
\]
so we have that 
\[
\mu _{o}=\exp _{(1,1,\delta _{f}(0))}, 
\]
as we claimed. $\blacksquare $

\bigskip

\section{\textbf{Lie product centered at an event}}

\bigskip

Now we desire also to see the evolution of any capitalized event as an
exponential map.

\bigskip

We define, at this aim, a new product induced by a capitalization $f$ and an
event $e_{0}$.

\bigskip

\textbf{Definition (Lie product induced by a capitalization factor and
centered at a point).}\emph{\ Let }$f$\emph{\ be a capitalization factor
over the real time line and let }$e_{0}$\emph{\ be any capitalized event }$%
(t_{0},h_{0},c_{0})$\emph{. Let }$e=(t,h,c)$\emph{\ and }$e^{\prime
}=(t^{\prime },h^{\prime },c^{\prime })$\emph{\ two capitalized financial
events, we define their }$(f,e_{0})$\emph{-\textbf{Lie product} }$%
[e|e^{\prime }]_{e_{0}}$\emph{\ to be the capitalized financial event} 
\[
\lbrack e|e^{\prime }]_{(f,e_{0})}=\left( t+(t^{\prime }-t_{0}),h+h^{\prime
}-h_{0},\frac{c}{f(h)}\frac{c^{\prime }}{f(h^{\prime })}\frac{c_{0}^{-1}}{%
f(-h_{0})}f(h+h^{\prime }-h_{0})\right) . 
\]

\bigskip

We will prove that the above product is indeed a Lie product, that is the
following theorem. We could follow the above proof, but we want to follow a
totally new and more interesting way, using also the above result.

\bigskip

We need first the concept of translation of a Lie group structure.

\bigskip

\section{\textbf{Translation of a Lie semigroup structure}}

\bigskip

\textbf{Theorem 4.} \emph{Let }$G$\emph{\ be a commutative Lie semigroup and
let }$e_{0}$\emph{\ one of its invertible elements, consider the product }$%
[.|.]_{e_{0}}$\emph{\ on the supporting set of the semigroup }$G$\emph{\
defined by} 
\[
\lbrack e|e^{\prime }]_{e_{0}}=ee^{\prime }e_{0}^{-1},
\]
\emph{for every }$e,e^{\prime }$\emph{\ in }$G$\emph{. Then, this new
product is a Lie product. Moreover, the neutral element of this new product
is the element }$e_{0}$\emph{\ and, denoted by }$G_{e_{0}}$\emph{\ the new
Lie semigroup, an element }$e$\emph{\ is invertible in }$G_{e_{0}}$\emph{\
if and only if it is invertible in the original }$G$\emph{\ and the inverse
of an invertible element }$e$\emph{\ in }$G_{e_{0}}$ \emph{is the element }$%
e^{-1}e_{0}^{2}$\emph{, where }$e^{-1}$ \emph{is the inverse of }$e$\emph{\
in }$G$\emph{.}

\emph{\bigskip }

\emph{Proof.} It is clear that the new operation is also associative and
commutative. Indeed, for instance the associativity is given by 
\begin{eqnarray*}
\lbrack [e|e^{\prime }]_{e_{0}}|e^{\prime \prime }]_{e_{0}} &=&[e|e^{\prime
}]_{e_{0}}e^{\prime \prime }e_{0}^{-1}= \\
&=&(ee^{\prime }e_{0}^{-1})e^{\prime \prime }e_{0}^{-1}= \\
&=&e(e^{\prime }e_{0}^{-1}e^{\prime \prime })e_{0}^{-1}= \\
&=&[e|[e^{\prime }|e^{\prime \prime }]_{e_{0}}]_{e_{0}},
\end{eqnarray*}
for every $e,e^{\prime },e^{\prime \prime }$ in $G$. For what concerns the
invertibility, we have 
\[
\lbrack e|e^{-1}e_{0}^{2}]_{e_{0}}=ee^{-1}e_{0}^{2}e_{0}^{-1}=e_{0}, 
\]
as we claimed. $\blacksquare $

\bigskip

\emph{The new Lie semigroup }$G_{e_{0}}$\emph{\ is called the translation of 
}$G$\emph{\ by }$e_{0}$\emph{\ and the new product the translation of the
product of }$G$\emph{\ by }$e_{0}$\emph{.}

\bigskip

\section{\textbf{Financial group translation}}

\bigskip

Well, we have exactly the following result.

\bigskip

\textbf{Theorem 5.}\emph{\ The financial product }$[.|.]_{(f,e_{o})}$\emph{\
is the translation by the event }$e_{0}$\emph{\ of the financial product }$%
[.|.]_{f}$\emph{, that is we have}
\[
\lbrack .|.]_{(f,e_{o})}=\tau _{e_{0}}[.|.]_{f}.
\]

\emph{\bigskip }

\emph{Proof.} Indeed, just recalling that 
\[
ee^{\prime }e^{\prime \prime }=\left( t+t^{\prime }+t^{\prime \prime
},h+h^{\prime }+h^{\prime \prime },\frac{cc^{\prime }c^{\prime \prime
}f(h+h^{\prime }+h^{\prime \prime })}{f(h)f(h^{\prime })f(h^{\prime \prime })%
}\right) , 
\]
we have 
\begin{eqnarray*}
\tau _{e_{0}}[.|.]_{f}(e,e^{\prime }) &=&[e|e^{\prime }]_{f}e_{0}^{-1}= \\
&=&(ee^{\prime })(-t_{0},-h_{0},c_{0}^{-1})= \\
&=&\left( t+t^{\prime }-t_{0},h+h^{\prime }-h_{0},\frac{cc^{\prime
}c_{0}^{-1}f(h+h^{\prime }-h_{0})}{f(h)f(h^{\prime })f(-h_{0})}\right) = \\
&=&[e|e^{\prime }]_{(f,e_{o})},
\end{eqnarray*}
as we claimed. $\blacksquare $

\bigskip

From which it immediately follows the claimed result about the product $%
[.|.]_{(f,e_{o})}$.

\bigskip

\textbf{Theorem 6.}\emph{\ The space }$\Bbb{R}^{3}$\emph{\ of capitalized
events is a }$C^{1}$\emph{\ Lie semigroup with respect to the standard
Euclidean differentiable structure on }$\Bbb{R}^{3}$\emph{\ and to the
centered product }$[.|.]_{(f,e_{o})}$\emph{, for every }$C^{1}$\emph{\
capitalization factor }$f$\emph{\ and any event }$e_{0}$\emph{. Moreover,
the neutral element of this product is the element }$e_{0}$\emph{.}

\bigskip

\section{\textbf{The evolution of an event}}

\bigskip

Consider again the evolution of an event $e_{0}$, that is the curve 
\[
\mu _{e_{0}}:\Bbb{R}\rightarrow \Bbb{R}^{3}:\mu
_{e_{0}}(t)=(t,h_{0}+t-t_{0},c_{0}f(-h_{0})f(h_{0}+t-t_{0})). 
\]

\bigskip

\textbf{Theorem 6.}\emph{\ Let }$o$\emph{\ be the unit event }$(0,0,1)$\emph{%
\ and let }$e_{0}$\emph{\ be any other event and let }$\mu _{e_{0}}$\emph{\
be its evolution. Then the evolution of }$e_{0}$\emph{\ is the double
translation of the evolution of the origin }$o$\emph{\ with respect to the
event }$e_{0}$\emph{\ itself and to the instant of time }$t_{0}$\emph{, that
is we have} 
\[
\mu _{e_{0}}=\tau _{e_{0}}\circ \tau _{t_{0}}(\mu _{o}).
\]

\bigskip

\emph{Proof.} We have
\begin{eqnarray*}
\mu _{o}(t-t_{0})e_{0} &=&(t-t_{0},t-t_{0},f(t-t_{0}))(t_{0},h_{0},c_{0})= \\
&=&(t,h_{0}+t-t_{0},f(t-t_{0})f(t_{0}-t)c_{0}f(-h_{0})f(t_{0}-t+h_{0}))= \\
&=&\mu _{e_{0}}(t),
\end{eqnarray*}
as we claimed. $\blacksquare $

\bigskip

Now we have the following result.

\bigskip

\textbf{Theorem 8.}\emph{\ Consider the translation of the usual addition on
the real line by a time }$t_{0}$\emph{\ and denote it by }$+_{t_{0}}$\emph{.
Then, the evolution }$\mu _{e_{0}}$\emph{\ of the event }$e_{0}$\emph{\ is
just an homomorphism of the group }$(\Bbb{R},+_{t_{0}})$\emph{\ into the
semigroup }$\Bbb{R}_{e_{0}}^{3}$\emph{.}

\emph{\bigskip }

\emph{First proof.} Indeed, we have 
\begin{eqnarray*}
\lbrack \mu _{e_{0}}(t)|\mu _{e_{0}}(t^{\prime })]_{(f,e_{o})} &=&\tau
_{e_{0}}[\mu _{u}(t-t_{0})e_{0}|\mu _{u}(t^{\prime }-t_{0})e_{0}]_{f}= \\
&=&[\mu _{u}(t-t_{0})e_{0}|\mu _{u}(t^{\prime }-t_{0})e_{0}]_{f}e_{0}^{-1}=
\\
&=&[\mu _{u}(t-t_{0})|\mu _{u}(t^{\prime }-t_{0})]_{f}e_{0}= \\
&=&\mu _{u}(t-t_{0}+t^{\prime }-t_{0})e_{0}= \\
&=&\mu _{e_{0}}(t-t_{0}+t^{\prime }),
\end{eqnarray*}
as we claimed. $\blacksquare $

\bigskip

\emph{Second direct proof (more complicated).} Indeed, we have 
\[
\mu _{e_{0}}(t+t^{\prime }-t_{0})=\left( t+t^{\prime
}-t_{0},h_{0}+t+t^{\prime }-2t_{0},c_{0}\frac{f(h_{0}+t+t^{\prime }-2t_{0})}{%
f(h_{0})}\right) , 
\]
and consider the capital evolution 
\[
M:\Bbb{R\rightarrow R}:M(t)=c_{0}\frac{f(h_{0}+t-t_{0})}{f(h_{0})}, 
\]
we so have, setting $h:=t-t_{0}$ and $h^{\prime }:=t^{\prime }-t_{0}$, that 
\begin{eqnarray*}
\lbrack \mu _{e_{0}}(t)|\mu _{e_{0}}(t^{\prime })]_{e_{0}}
&=&[(t,h_{0}+h,M(t))|(t^{\prime },h_{0}+h^{\prime },M(t^{\prime
}))]_{(f,e_{0})}= \\
&=&\left( t+t^{\prime }-t_{0},h_{0}+h+h^{\prime },M(t)M(t^{\prime })\frac{%
f(h_{0}+h+h^{\prime })f(h_{0})}{c_{0}f(h_{0}+h)f(h_{0}+h^{\prime })}\right) =
\\
&=&(t+t^{\prime }-t_{0},h_{0}+t+t^{\prime
}-2t_{0},c_{0}f(-h_{0})f(h_{0}+t+t^{\prime }-2t_{0})),
\end{eqnarray*}
for every $t$ and $t^{\prime }$ on the time line. $\blacksquare $

\bigskip

\textbf{Theorem 9.}\emph{\ The evolution of the origin }$e_{0}$\emph{\ of
the Lie semigroup }$\Bbb{R}_{e_{0}}^{3}$\emph{\ is the exponential map of
the translated Lie group }$\Bbb{R}_{t_{0}}$ \emph{into }$\Bbb{R}_{e_{0}}^{3}$%
\emph{, with respect to the tangent vector }$(e_{0},(1,1,c_{0}\delta
(h_{0})))$\emph{, tangent vector to the Lie semigroup }$\Bbb{R}_{e_{0}}^{3}$%
\emph{\ at its own origin }$e_{0}$\emph{.}

\emph{\bigskip }

\emph{Proof.} The tangent vector at $t$ of $\mu _{e_{0}}$ is 
\[
\mu _{e_{0}}^{\prime }(t)=(1,1,c_{0}f(-h_{0})f^{\prime }(h_{0}+t-t_{0})),
\]
so that we have 
\[
\mu _{e}^{\prime }(t_{0})=(1,1,c_{0}f^{\prime
}(h_{0})/f(h_{0}))=(1,1,c_{0}\delta _{f}(h_{0})),
\]
where 
\[
\delta _{f}:\Bbb{R}\rightarrow \Bbb{R}:h\mapsto f^{\prime }(h)/f(h),
\]
is the force of interest of the capitalization factor $f$. Now, there is
only one $1$-parameter group in $\Bbb{R}_{e_{0}}^{3}$ having $v$ as a
tangent vector at the origin $t_{0}$ of the group $(\Bbb{R},+_{t_{0}})$, and
it is the exponential map 
\[
\exp _{v}:\Bbb{R}_{t_{0}}\rightarrow \Bbb{R}_{e_{0}}^{3},
\]
so we have 
\[
\mu _{e_{0}}=\exp _{(1,1,c_{0}\delta _{f}(h_{0}))},
\]
as we claimed. $\blacksquare $

\bigskip

\section{\textbf{General exponential map on a commutative Lie group}}

\bigskip

Let $G$ be a Lie group we know that for every tangent vector $v$ at the
origin (that is for every element $v$ of its associated Lie algebra) there
is a unique homomorphism $\mu :\Bbb{R}\rightarrow G$ such that 
\[
d_{0}\mu (h)=hv, 
\]
for every $h$ in $\Bbb{R}$.

\bigskip

Now, for every pair $(t_{0},e_{0})$ in the product $\Bbb{R}\times G$,
consider the double translation 
\[
\mu _{(t_{0},e_{0})}:t\mapsto \mu (t-t_{0})e_{0}, 
\]
it is clear that $\mu _{(t_{0},e_{0})}$ is a homomorphism from the
translated group $\Bbb{R}_{t_{0}}$ into the translated $G_{e_{0}}$.
Moreover, we have 
\[
d_{t_{0}}\mu _{(t_{0},e_{0})}(h)=d_{0}\mu (h) 
\]

\bigskip

\textbf{Definition (the general exponential).}\emph{\ We call exponential
map of the Lie subgroup }$G$\emph{\ at }$(t_{0},e_{0})$\emph{\ relative to
the tangent vector }$v$\emph{\ in }$T_{e_{0}}(G)$\emph{\ the unique
homomorphism }$\mu $\emph{\ from translated group }$\Bbb{R}_{t_{0}}$\emph{\
into the translated }$G_{e_{0}}$\emph{\ such that} 
\[
d_{t_{0}}\mu (h)=hv, 
\]
\emph{for every }$h$\emph{\ in }$\Bbb{R}$\emph{\ (note that the application }%
$d_{t_{0}}\mu $\emph{\ goes from }$T_{t_{0}}(\Bbb{R})$\emph{\ into }$%
T_{e_{0}}(G_{e_{0}})$\emph{\ since }$\mu (t_{0})$\emph{\ is }$e_{0}$\emph{).
We denote this }$\mu $\emph{\ by } 
\[
\exp _{(t_{0},e_{0}),v}. 
\]

\bigskip

\textbf{Remark.} Note that the tangent space $T_{e_{0}}(G_{e_{0}})$ is the
tangent space $T_{e_{0}}(G)$, indeed a derivation $v$ in $%
T_{e_{0}}(G_{e_{0}})$ is a functional defined on $C_{e_{0}}^{1}(G,\Bbb{R})$
and so it is a derivation in $T_{e_{0}}(G)$.

\bigskip

With the above definition, we can say that

\begin{itemize}
\item  \emph{the evolution of a capitalized event }$e_{0}$\emph{\ is the
exponential of the Lie semigroup }$\Bbb{R}_{f}^{3}$\emph{\ at the pair }$%
(t_{0},e_{0})$\emph{\ with respect to the tangent vector }$%
(e_{0},(1,1,c_{0}\delta _{f}(h_{0}))$\emph{.}
\end{itemize}

\bigskip

\textbf{David Carf\`{i}}

\emph{Faculty of Economics}

\emph{University of Messina}

\emph{davidcarfi71@yahoo.it}

\end{document}